\documentclass[12pt]{article}
\usepackage{amssymb}
\usepackage{amsmath}
\usepackage{longtable}

\newcommand{\be}{\begin{equation}}
\newcommand{\ee}{\end{equation}}
\newcommand{\bea}{\begin{eqnarray*}}
\newcommand{\eea}{\end{eqnarray*}}
\newcommand{\ba}{\begin{array}}
\newcommand{\ea}{\end{array}}
\newcommand{\bi}{\begin{itemize}}
\newcommand{\ei}{\end{itemize}}
\newcommand{\bc}{\begin{center}}
\newcommand{\ec}{\end{center}}
\newcommand{\bfr}{\begin{flushright}}
\newcommand{\efr}{\end{flushright}}
\newcommand{\f}{\frac}
\newcommand{\ov}{\overline}

\newcommand{\eps}{\varepsilon}

\newcommand{\ds}{\displaystyle}

\newcommand{\q}{\quad}

\begin{document}

\title{On bi-unitary harmonic numbers}
\author{J\'ozsef S\'andor\\
Babe\c{s}-Bolyai University\\
Faculty of Mathematics and Informatics\\
Str. Kog\u{a}lniceanu 1\\
400084 Cluj-Napoca, Romania}
\date{}
\maketitle

\begin{abstract}
The aim of this paper is twofold. First we give a short survey of the
existing results on various notions of harmonic numbers; and then we make
a preliminary study of bi-unitary harmonic numbers.
\end{abstract}

\section{Introduction}
In 1948 O. Ore \cite{9} considered numbers $n$ whose divisors
have integral harmonic mean
$$H(n)=\f{r}{\ds\f{1}{d_1}+\dots +\ds\f{1}{d_r}},$$
where $1=d_1<d_2<\dots <d_r=n$
are all the divisors of $n$.
Since
$$\sum\f{1}{d_r}=\f{1}{n}\sum \f{n}{d_r}=\f{1}{n}\sigma (n),$$
and $r=d(n)$ (with $\sigma (n)$ and $d(n)$ denoting the sum, resp.
number of divisors of $n$), clearly
$$H(n)=\f{nd(n)}{\sigma (n)},$$
so $H(n)$ is an integer iff
$$\sigma (n)|nd(n)
\eqno(1)$$

C. Pomerance \cite{10} called a number $n$ with property (1),
a {\bf harmonic number}.
Ore proved that if $n$ is perfect (i.e. $\sigma (n)=2n$),
then it is harmonic.
Indeed, if $n$ is perfect, then $2|d(n)$ is always true,
i.e. $n$ is not a perfect square.

Ore proved also that if $n$ is harmonic, then
$\omega (n)\ge 2$, and Pomerance showed that the only harmonic
numbers with two distinct prime factors are the even perfect numbers.
In 1963 M. V. Subbarao \cite{11} called the number $n$ a
{\bf balanced number} if
$\ds\f{\sigma (n)}{d(n)}=\ds\f{n}{2}$,
and proved that $n=6$ is the single balanced number.
Now, remark that a balanced number satisfies $H(n)=2$,
so it is a particular harmonic number.
M. Garcia \cite{4} extended the list of harmonic numbers to include
all 45 which are $<10^7$, and found more than 200 larger ones.
The least one, apart from 1 and the perfect numbers, is 140.
All 130 harmonic numbers up to $2\cdot 10^9$ are listed by
G. L. Cohen \cite{1}; and R. M. Sorli (see Cohen and Sorli \cite{2})
has continued the list to $10^{10}$.
Ore conjectured that every harmonic number is even, but this
probably is very difficult.
Indeed, this result, if true, would imply that there are no odd
perfect numbers.
See also W. H. Mills \cite{7}, who proved that if there exists
an odd harmonic number $n$, then $n$ has  a prime-power factor
greater than $10^7$.

In 1998 G. L. Cohen and M. Deng \cite{3} have introduced
a generalization of harmonic numbers.
Let $k\ge 1$, integer and let $\sigma _k(n)$ be the sum of $k$th
powers of divisors of $n$.
Then $n$ is called {\bf $k$-harmonic}, if
$$\sigma _k(n)|n^k d(n)
\eqno(2)$$

They proved that for $k>1$ there is no $k$-harmonic number
in the range $1<n\le 10^{10}$.

A divisor $d$ of $n$ is called {\bf unitary divisor}, if
$\left(d,\f{n}{d}\right)=1$.
Let $\sigma ^*(n)$ and $d^*(n)$ denote the sum, resp. number- of
unitary divisors of $n$.
M. V. Subbarao and L. J. Warren \cite{12} introduced the
{\bf unitary perfect} numbers $n$ satisfying
$\sigma ^*(n)=2n$.
They found the first four unitary perfect numbers, while the fifth
one was discovered by Ch. Wall \cite{15}.
A number $n$ is called {\bf unitary harmonic} if
$$\sigma ^*(n)|nd^*(n),
\eqno(3)$$
concept introduced by K. Nageswara Rao \cite{8}, who showed that if $n$
is unitary perfect, then it is also unitary harmonic.
P. Hagis and G. Lord \cite{5} proved that if $H^*(x)$
is the counting function of these numbers, then for $\eps >0$
and large $x$ one has
$$H^*(x)<2.2x^{1/2}\cdot 2^{(1+\eps )\log x/\log\log x}$$

The same result was obtained in 1957 by H.-J. Kanold [8]
for the counting function of harmonic numbers.

Wall [19] showed that there are 23 unitary harmonic numbers $n$
with $\omega (n)\le 4$, and claimed that there are 43 unitary harmonic
numbers $n\le 10^6$.
However, Hagis and Lord \cite{5} have shown this with 45 in place of 43.

Recently, T. Goto and S. Shibata [5] have determined all harmonic
numbers satisfying $H(n)\le 300$.
According to the Referee, Goto and K. Okeya have extended the study
to $H(n)\le 1200$.

For {\bf infinitary harmonic} numbers, related to the concept
of an "infinitary divisor", see Hagis and Cohen [7].

For many results involving these topics, see also Chapter I
of the author's recent book [13].

\section{Bi-unitary divisors and bi-unitary harmonic numbers}
A divisor $d$ of $n$ is called a {\bf bi-unitary divisor} if the
greatest common unitary divisor of $d$ and $\ds\f{n}{d}$ is 1.
Let $\sigma ^{**}(n)$ be the sum of bi-unitary divisors of $n$.
Wall \cite{14} called a number $n$ {\bf bi-unitary perfect},
if $\sigma ^{**}(n)=2n$,
and proved that there are only three bi-unitary perfect numbers,
namely 6, 60 and 90.
It is not difficult to verify that $\sigma ^{**}$ is multiplicative,
and for prime powers $p^\alpha $,
$$\ba{llll}
\sigma ^{**}(p^\alpha )=\sigma (p^\alpha )=
\ds\f{p^{\alpha +1}-1}{p-1},
& \mbox{if} & \alpha & \mbox{is odd}
\\
\sigma ^{**}(p^\alpha )=\ds\f{p^{\alpha +1}-1}{p-1}-p^{\alpha /2},
& \mbox{if} & \alpha & \mbox{is even}
\ea$$

Let $d^{**}(n)$ be the number of bi-unitary divisors of $n$.
Then it is also known (see D. Suryanarana \cite{13}) that if
$n=p_1^{a_1}\dots p_r^{a_r}$
is the prime factorization of $n>1$, then
$$d^{**}(n)=\left(\prod_{a_i=even}a_i\right)
\left(\prod_{a_i=odd}(a_i+1)\right)
\eqno(5)$$

We now introduce the main notion and results of this paper.

{\bf Definition.}
The number $n$ is called {\bf bi-unitary harmonic}, if
$$\sigma ^{**}(n)|nd^{**}(n)
\eqno(6)$$

{\bf Theorem 1.}
{\it Let $k\ge 1$ be an integer and suppose that $n$
is {\bf bi-unitary $k$-perfect}, i.e.
$\sigma ^{**}(n)=kn$.
Then $n$ is bi-unitary harmonic iff
$$k|d^{**}(n)
\eqno(7)$$

Particularly, 6, 60, 90 are bi-unitary harmonic numbers.
}

{\bf Proof.}
(7) is a consequence of Definition (6) and the definition of
bi-unitary $k$-perfect numbers.
Remark that for $k=2$, by relation (5), (7) is always true.
By Wall's result on bi-unitary perfect numbers it follows that 6, 60, 90
are also bi-unitary harmonic numbers.

Let $\omega (n)$ denote the number of distinct prime factors of $n$.

{\bf Corollary.}
{\it If $\omega (n)\ge 2$ and $n$ is bi-unitary
4-perfect, then $n$ is bi-unitary harmonic number.
}

{\bf Proof.}
Remark that for $\omega (n)\ge 2$, by (5),
$4|d^{**}(n)$ so by (7) the result follows.

{\bf Theorem 2.}
{\it Let $n=p_1^{a_1}\dots p_r^{a_r}>1$
be the prime factorization of $n$, and suppose that all $a_i$
$(i=\ov{1,r})$ are odd.
Then $n$ is bi-unitary harmonic iff it is harmonic.
}

{\bf Proof.}
Since all $a_i$ $(i=\ov{1,r})$ are odd, by (5),
$$d^{**}(n)=\prod_{a_i=odd}(a_i+1)=d(n),$$
and by (4),
$$\sigma ^{**}(n)=\prod_{\alpha \, odd}\sigma ^{**}(p^\alpha )
=\prod_{\alpha \, odd}\sigma (p^\alpha )=\sigma (n).$$

Thus (6) is true if and only if (1) is true.

{\bf Corollary.}
1) Besides 1, the only squarefree bi-unitary harmonic number is 6.
This follows by a result of Ore \cite{9} on harmonic numbers.

2) If $n$ is odd bi-unitary harmonic, with all $a_i$ odd,
then $n$ has a component exceeding $10^7$.
If $n$ is even, then $\omega (n)\ge 3$.

This follows by the result of Mills stated in the Introduction,
as well as by the fact that if
$\omega (n)=2$, then $n$ being harmonic, it must be a perfect number.
All even perfect number has the form $2^{2k}p$, where $p$
is an odd prime.
Since $2k=$ even, this leads to a contradiction.

{\bf Remark.}
By Theorem 2, new bi-unitary harmonic numbers can be found.
For example,
$n=2^1\cdot 3^3\cdot 5=270$,
$n=2^5\cdot 3\cdot 7=672$,
$n=2^5\cdot 3^3\cdot 5\cdot 7=30240$,
$n=2^3\cdot 3^3\cdot 5^3\cdot 7\cdot 13=2457000$,
see \cite{2} for a list of harmonic seeds and harmonic numbers.

A computer program, on the other hand, may be applied for a search
of bi-unitary harmonic numbers.
For example, there are 50 such numbers $n\le 10^6$, but the search
could be extended to $10^9$ (see the Table at the end of the paper), etc.

{\bf Theorem 3.}
{\it Let $n=p_1^{a_1}\dots p_r^{a_r}>1$
be the prime factorization of $n$.
If all $a_i$ $(i=\ov{1,r})$ are even, and $n$ is bi-unitary harmonic,
then
$$\omega (n)\ge 2
\eqno(8)$$
}

{\bf Proof.}
If $n=p^{2a}$ is a prime power, with even exponent $2a$,
then by (5), (6) is equivalent to
$(1+p+p^2+\dots +p^{a-1}+p^{a+1}+\dots +p^{2a})|p^{2a}(2a)$.
Clearly $p^{2a}$ is relatively prime to
$1+p+\dots +p^{a-1}+p^{a+1}+\dots +p^{2a}$,
so we must have
$$(1+p+p^2+\dots +p^{a-1}+p^{a+1}+\dots +p^{2a})|(2a)$$

But this is impossible, since the first term contains a number of
$2a$ terms, each (excepting 1) greater than 1, so
$1+p+\dots +p^{a-1}+p^{a+1}+\dots +p^{2a}>2a$.
Therefore (8) follows.

{\bf Corollary.}
{\it If $n$ is bi-unitary harmonic number, then $\omega (n)\ge 2$.
}

Indeed, by Theorem 3, $n$ cannot be of the form $p^{2a}$.
On the other hand, if
$n=p^{2a+1}$ ($p$ prime), then it is harmonic, contradicting Ore's
result that $\omega (n)\ge 2$.

If there are odd, as well as even exponents, the following particular
result holds true:

{\bf Theorem 4.}
{\it There are no bi-unitary harmonic numbers of the form
$p^3\cdot q^2$ ($p,q$ distinct primes).
}

{\bf Proof.}
If $n=p^{2a+1}q^{2b}$, then
$$d^{**}(n)=2b(2a+2)=4b(a+1),$$
$$\sigma ^{**}(n)=(1+p+\dots +p^{2a+1})(1+q+\dots +q^{b-1}+q^{b+1}
+\dots +q^{2b}),$$
so $n$ is bi-unitary harmonic iff
$$(1+p+\dots +p^{2a+1})(1+q+\dots +q^{b-1}+q^{b+1}+\dots +q^{2b})|
4b(a+1)p^{2a+1}q^{2b}
\eqno(9)$$

For $a=1$, $b=1$ this becomes
$$(1+p+p^2+p^3)(1+q^2)|8p^3q^2
\eqno(10)$$

Since
$(1+p+p^2+p^3,p^3)=1$ and $(1+q^2,q^2)=1$,
it follows that
$$(1+p+p^2+p^3)|8q^2\mbox{ and } (1+q^2)|8p^3
\eqno(11)$$

If $p=2$, it follows
$(1+q^2)|64$, which is impossible for all $q$.
Similarly, if $q=2$, then $5|8p^3$, so $p=5$, and then the first relation
of (11) is impossible.
Thus, $p,q\ge 3$.
Remark that
$1+p+p^2+p^3=(1+p)(1+p^2)$ and that (11) implies that
$(1+p)$ and $(1+p^2)$ can have only two distinct prime factors,
namely: 2 and $q$.
Let $1+p=kq$, i.e. $p=kq-1$.
Then $p^2+1=k^2q^2-2kq+2$ is divisible by $q$ only if $q=2$.
But this is impossible.
Clearly $k$ is even, $k=2s$, so $p^2+1=2(2s^2q^2-2sq+1)$
which cannot be a power of 2, since
$2s^2q^2-2sq+1$ is odd.
Thus $(1+p)(1+p^2)$ can have also other prime factors than 2 and $q$,
contradicting (11).

{\bf Remark.}
It can be proved similarly that there are no bi-unitary harmonic
numbers of the form $pq^4$
or $p^3\cdot q^4$, and that the only one of the form
$pq^2$ is $5\cdot 3^2$.

Finally, we state the following result:

{\bf Theorem 5.}
{\it Let $n$ be of the form
$$n=\ds\prod_{i=1}^r p_i^{2a_i+1}\prod_{j=1}^s q_j^{2b_j},$$
and let
$$n_1=\prod_{i=1}^r p_i^{2a_i+1}\prod_{j=1}^s q_j^{b_j-1}
\q\mbox{and}\q
n_2=\prod_{j=1}^s q_j^{b_j}$$
($p_i$ and $q_j$ are distinct primes;
$a_i,b_j$ positive integers).
Suppose that $n_1$ is a harmonic number, while $n_2$ a unitary harmonic
number.
Then $n=n_1n_2$ is a bi-unitary harmonic number.
}

{\bf Proof.}
This follows from the fact that for the numbers $n$ given above
one has the identity
$$H^{**}(n)=H(n_1)H^*(n_2),
\eqno(12)$$
where $H,H^*,H^{**}$ are the corresponding harmonic means
(e.g. $H^{**}(n)=nd^{**}(n)/\sigma ^{**}(n)$).
Identity (12) can be proved by using the definitions and
the results (e.g. relation (5)) for the above functions.

{\bf Final remarks.}
By examining the Table with all bi-unitary harmonic numbers
up to $10^9$, we can remark that there are in total 211 such numbers
in this range.
The first 12 bi-unitary harmonic numbers are all harmonic,
or unitary harmonic; the first number without this property is
$n=9072$.
There are only 5 bi-unitary harmonic numbers up to $10^9$
which are powerful, namely
$n=3307500,\ 9922500,\ 23152500,\ 138915000,\ 555660000$.
From these, only
$n=9922500=2^2\cdot 3^4\cdot 5^4\cdot 7^2$
is a perfect square.
It is interesting to note that, the existence of perfect squares
in the set of harmonic or unitary harmonic numbers, is an open
question up to now.

\section{Related numbers}
As we have seen, the harmonic means of divisors, unitary
divisors, and bi-unitary divisors are given explicitly by
$$H(n)=\f{nd(n)}{\sigma (n)},\q
H^*(n)=\f{nd^*(n)}{\sigma ^*(n)},\q
H^{**}(n)=\f{nd^{**}(n)}{\sigma ^{**}(n)}.$$

In what follows, the harmonic, unitary harmonic, resp.
bi-unitary harmonic numbers will be called simply as $H,H^*$,
resp. $H^{**}$-numbers.
This will be motivated also by the introduction of the following
six new fractions and related numbers:
$$H_1(n)=\f{nd(n)}{\sigma ^*(n)},\q
H_2(n)=\f{nd^*(n)}{\sigma (n)},\q
H_3(n)=\f{nd(n)}{\sigma ^{**}(n)},$$
$$H_4(n)=\f{nd^{**}(n)}{\sigma (n)},\q
H_5(n)=\f{nd^*(n)}{\sigma ^{**}(n)},\q
H_6(n)=\f{nd^{**}(n)}{\sigma ^*(n)}.$$

When $H_1(n)$ is an integer, we will say that $n$
{\bf is a $H_1$ number}, etc.

By remarking that $d^{**}(n)$ is always divisible by $d^*(n)$,
and that if $n$ has the form
$n=p_1^{\eps _1}\dots p_r^{\eps _r}$, where
$\eps _i\in\{1,2\}$ $(i=\ov{1,r})$,
$p_i$ distinct primes, then
$d^{**}(n)=d^*(n)$, $\sigma ^{**}(n)=\sigma ^*(n)$,
we can state the following result:

{\bf Theorem 6.}
{\it In all cases,
$H^{**}(n)=H_5(n)k_1(n)$,
$H_4(n)=H_2(n)k_2(n)$,
$H_6(n)=H^*(n)k_3(n)$,
where
$k_1(n),k_2(n),k_3(n)$
are integers.
If $n$ has the form
$n=p_1^{\eps _1}\dots p_r^{\eps _r}$, then
$H_3(n)=H_1(n)$, $H_5(n)=H^*(n)$, $H_6(n)=H^{**}(n)=H^*(n)$,
$H_4(n)=H_2(n)$.

If $n$ has the form
$p_1^{a_1}\dots p_r^{a_r}$
with all $a_i$ odd $(i=\ov{1,r})$, then
$H_5(n)=H_2(n)$, $H_6(n)=H_1(n)$, $H_3(n)=H_4(n)=H(n)$.
}

{\bf Corollary.}
In all cases, a $H_2$-number is also a $H_4$-number;
a $H^*$-number is also a $H_6$-number;
a $H_5$-number is also a $H^{**}$-number.
If $n=p_1^{\eps _1}\dots p_r^{\eps _r}$, then $n$ is
a $H^{**}$-number iff it is a $H^*$, and a $H_6$-number;
$n$ is $H_5$-number, iff it is a $H^*$-number, and
$n=H_2$-number iff $n=H_4$-number.
If all $a_i$ are odd, then the notions of $H_2$ and $H_5$-numbers;
$H,H_3,H_4$-numbers; resp. $H_1$ and $H_6$-numbers coincide.

{\bf Remark.}
Since in Wall \cite{16} there are stated all $H^*$-numbers with
$\omega (n)\le 4$, we can say from the Table 1 of that paper,
that the only $H^{**}$ (i.e. bi-unitary harmonic)-number of the form
$p^2q$ is
$3^2\cdot 5=45$;
- of the form $p^2qr$ are $2^2\cdot 3\cdot 5=60$ and
$3^2\cdot 2\cdot 5=90$;
- of the form $p^2q^2r$ is $5^2\cdot 7^2\cdot 13=15925$;
- of the form $p^2qrs$ are $2^2\cdot 3\cdot 5\cdot 7=420$,
$3^2\cdot 2\cdot 5\cdot 7=630$;
- of the form $p^2q^2rs$ is $2^2\cdot 5^2\cdot 7\cdot 13=9100$;
- of the form $p^2q^2r^2s$ is $3^2\cdot 5^2\cdot 13^2\cdot 17=646425$.
(Here $p,q,r,s$ are distinct primes).
These complement some results of Theorem 4.

Clearly, the deeper study of all of the above numbers cannot be done
in this paper (but there are some results under preparation).
We state only the following result:

{\bf Theorem 7.}
{\it If $n$ is a perfect number, then $n$ is a $H_2$ and $H_4$-number, too.
If $n>1$ is a $H_2$-number, then it cannot be a geometric number
(i.e. a perfect square).
}

{\bf Proof.}
Let $\sigma (n)=2n$.
Then, as $d^*(n)=2^{\omega (n)}$, and
$2|d^{**}(n)$, clearly $H_2(n)$ and $H_4(n)$ will be integers.
It is well known that $\sigma (n)$ is odd if $n$ is a perfect square
(i.e., $n=m^2$).
Then, if $H_2(n)$ is an integer, then clearly $\sigma (n)$
divides $n$, and this is possible only for $n=1$.

{\bf Remarks.}
1) The similar problem in the case of $H$-numbers, i.e.
if they are or not geometric, is a difficult open question
(see e.g. \cite{2}).

2) A number $n>1$ is called {\bf friendly number}
(or Duffinian number), see e.g. \cite{19}, if
$(n,\sigma (n))=1$.
Clearly, if $n$ is friendly, then $n$ cannot be $H_2$- or $H_4$-number.
Indeed, in this case one must have
$\sigma (n)|d^*(n)$, or
$\sigma (n)|d^{**}(n)$, but this is impossible for $n>1$, since
$\sigma (n)>d^{**}(n)\ge d^*(n)$ for $n>1$.
A similar result holds true for the $H$-numbers.

{\bf Acknowledgements.}
The author is indebted to the Referee
as well as to Professor G. L. Cohen,
for pointing out new references;
for corrections in a former version of the paper;
and for many suggestions which considerably improved the presentation
of the paper.
He also thanks Professor R. M. Sorli for providing him a list of
bi-unitary harmonic numbers up to $10^9$ (see the attached Table).

\bigskip

\noindent
AMS Subject classification: 11A25, 11A99, 11N37.

\newpage

\bc
All bi-unitary harmonic numbers less than $10^9$
\ec

\bc
\begin{longtable}{|l|l|}\hline
$n$ & $H^{**}(n)$\\ \hline \endhead
1 & 1\\
$6=2\cdot 3$ & 2\\
$45=3^2\cdot 5$ & 3\\
$60=2^2\cdot 3\cdot 5$ & 4\\
$90=2\cdot 3^2\cdot 5$ & 4\\
$270=2\cdot 3^3\cdot 5$ & 6\\
$420=2^2\cdot 3\cdot 5\cdot 7$ & 7\\
$630=2\cdot 3^2\cdot 5\cdot 7$ & 7\\
$672=2^5\cdot 3\cdot 7$ & 8\\
$2970=2\cdot 3^3\cdot 5\cdot 11$ & 11\\
$5460=2^2\cdot 3\cdot 5\cdot 7\cdot 13$ & 13\\
$8190=2\cdot 3^2\cdot 5\cdot 7\cdot 13$ & 13\\
$9072^*=2^4\cdot 3^4\cdot 7$ & 12\\
$9100=2^2\cdot 5^2\cdot 7$ & 10\\
$10080=2^5\cdot 3^2\cdot 5\cdot 7$ & 16\\
$15925=5^2\cdot 7^2\cdot 13$ & 7\\
$22680=2^3\cdot 3^4\cdot 5\cdot 7$ & 18\\
$22848=2^6\cdot 3\cdot 7\cdot 17$ & 16\\
$27300=2^2\cdot 3\cdot 5^2\cdot 7\cdot 13$ & 15\\
$30240=2^5\cdot 3^3\cdot 5\cdot 7$ & 24\\
$40950=2\cdot 3^2\cdot 5^2\cdot 7\cdot 13$ & 15\\
$45360=2^4\cdot 3^4\cdot 5\cdot 7$ & 20\\
$54600=2^3\cdot 3\cdot 5^2\cdot 7\cdot 13$ & 20\\
$81900=2^2\cdot 3^2\cdot 5^2\cdot 7\cdot 13$ & 18\\
$95550=2\cdot 3\cdot 5^2\cdot 7^2\cdot 13$ & 14\\
$99792=2^4\cdot 3^4\cdot 7\cdot 11$ & 22\\
$136500=2^2\cdot 3\cdot 5^3\cdot 7\cdot 13$ & 25\\
$163800=2^3\cdot 3^2\cdot 5^2\cdot 7\cdot 13$ & 24\\
$172900=2^2\cdot 5^2\cdot 7\cdot 13\cdot 19$ & 19\\
$204750=2\cdot 3^2\cdot 5^3\cdot 7\cdot 13$ & 25\\
$208656=2^4\cdot 3^4\cdot 7\cdot 23$ & 23\\
$245700=2^2\cdot 3^3\cdot 5^2\cdot 7\cdot 13$ & 27\\
$249480=2^3\cdot 3^4\cdot 5\cdot 7\cdot 11$ & 33\\
$312480=2^5\cdot 3^2\cdot 5\cdot 7\cdot 31$ & 31\\
$332640=2^5\cdot 3^3\cdot 5\cdot 7\cdot 11$ & 44\\
$342720=2^6\cdot 3^2\cdot 5\cdot 7\cdot 17$ & 32\\
$385560=2^3\cdot 3^4\cdot 5\cdot 7\cdot 17$ & 34\\
$409500=2^2\cdot 3^2\cdot 5^3\cdot 7\cdot 13$ & 30\\
$472500=2^2\cdot 3^3\cdot 5^4\cdot 7$ & 25\\
$491400=2^3\cdot 3^3\cdot 5^2\cdot 7\cdot 13$ & 36\\
$646425=3^2\cdot 3^2\cdot 13^2\cdot 17$ & 13\\
$695520=2^5\cdot 3^3\cdot 5\cdot 7\cdot 23$ & 46\\
$708288=2^6\cdot 3\cdot 7\cdot 17\cdot 31$ & 31\\
$716625=3^2\cdot 5^3\cdot 7^2\cdot 13$ & 21\\
$791700=2^2\cdot 3\cdot 5^2\cdot 7\cdot 13\cdot 29$ & 29\\
$819000=2^3\cdot 3^2\cdot 5^3\cdot 7\cdot 13$ & 40\\
$861840=2^4\cdot 3^4\cdot 5\cdot 7\cdot 19$ & 38\\
$900900=2^2\cdot 3^2\cdot 5^2\cdot 7\cdot 11\cdot 13$ & 33\\
$955500=2^2\cdot 3\cdot 5^3\cdot 7^2\cdot 13$ & 28\\
$982800=2^4\cdot 3^3\cdot 5^2\cdot 7\cdot 13$ & 40\\
$1028160=2^6\cdot 3^3\cdot 5\cdot 7\cdot 17$ & 48\\
$1037400=2^3\cdot 3\cdot 5^2\cdot 7\cdot 13\cdot 19$ & 38\\
$1187550=2\cdot 3^2\cdot 5^2\cdot 7\cdot 13\cdot 29$ & 29\\
$1228500=2^2\cdot 3^3\cdot 5^3\cdot 7\cdot 13$ & 45\\
$1392300=2^2\cdot 3^2\cdot 5^2\cdot 7\cdot 13\cdot 17$ & 34\\
$1421280=2^5\cdot 3^3\cdot 5\cdot 7\cdot 47$ & 47\\
$1433250=2\cdot 3^2\cdot 5^3\cdot 7^2\cdot 13$ & 28\\
$1528800=2^5\cdot 3\cdot 5^2\cdot 7^2\cdot 13$ & 32\\
$1571328=2^9\cdot 3^2\cdot 11\cdot 31$ & 32\\
$1801800=2^3\cdot 3^2\cdot 5^2\cdot 11\cdot 13$ & 44\\
$2457000=2^3\cdot  3^3\cdot 5^3\cdot 7\cdot 13$ & 60\\
$2579850=2\cdot 3^4\cdot 5^2\cdot 7^2\cdot 13$ & 27\\
$2888704=2^{10}\cdot 7\cdot 13\cdot 31$ & 32\\
$3307500^*=2^2\cdot 3^3\cdot 5^4\cdot 7^2$ & 28\\
$3767400=2^3\cdot 3^2\cdot 5^2\cdot 7\cdot 13\cdot 23$ & 46\\
$3878550=2\cdot 3^3\cdot 5^2\cdot 13^2\cdot 17$ & 26\\
$4176900=2^2\cdot 3^3\cdot 5^2\cdot 7\cdot 13\cdot 17$ & 51\\
$4291056=2^4\cdot 3^4\cdot 7\cdot 11\cdot 43$ & 43\\
$4299750=2\cdot 3^3\cdot 5^3\cdot 7^2\cdot 13$ & 42\\
$4504500=2^2\cdot 3^2\cdot 5^3\cdot 7\cdot 11\cdot 13$ & 55\\
$4713984=2^9\cdot 3^3\cdot 11\cdot 31$ & 48\\
$4961250=2\cdot 3^4\cdot 5^4\cdot 7^2$ & 25\\
$5405400=2^3\cdot 3^3\cdot 5^2\cdot 7\cdot 11\cdot 13$ & 66\\
$6168960=2^7\cdot 3^4\cdot 5\cdot 7\cdot 17$ & 64\\
$6397300=2^2\cdot 5^2\cdot 7\cdot 13\cdot 19\cdot 37$ & 37\\
$6688500=2^2\cdot 3\cdot 5^3\cdot 7^3\cdot 13$ & 49\\
$7698600=2^3\cdot 3^2\cdot 5^2\cdot 7\cdot 13\cdot 47$ & 47\\
$7780500=2^2\cdot 3^2\cdot 5^3\cdot 7\cdot 13\cdot 19$ & 57\\
$7983360=2^8\cdot 3^4\cdot 5\cdot 7\cdot 11$ & 64\\
$8353800=2^3\cdot 3^3\cdot 5^2\cdot 7\cdot 13\cdot 17$ & 68\\
$8666112=2^{10}\cdot 3\cdot 7\cdot 13\cdot 31$ & 48\\
$9922500^*=2^2\cdot 3^4\cdot 5^4\cdot 7^2$ & 30\\
$10032750=2\cdot 3^2\cdot 5^3\cdot 7^3\cdot 13$ & 49\\
$10624320=2^6\cdot 3^2\cdot 5\cdot 7\cdot 17\cdot 31$ & 62\\
$10701600=2^5\cdot 3\cdot 5^2\cdot 7^3\cdot 13$ & 56\\
$10999296=2^9\cdot 3^2\cdot 7\cdot 11\cdot 31$ & 56\\
$11302200=2^3\cdot  3^3\cdot 5^2\cdot 7\cdot 13\cdot 23$ & 69\\
$11309760=2^6\cdot  3^3\cdot 5\cdot 7\cdot 11\cdot 17$ & 88\\
$11875500=2^2\cdot  3^2\cdot 5^3\cdot 13\cdot 29$ & 58\\
$12899250=2\cdot 3^4\cdot 5^3\cdot 7^2\cdot 13$ & 45\\
$13022100=2^2\cdot 3^3\cdot 5^2\cdot 7\cdot 13\cdot 53$ & 53\\
$14303520=2^5\cdot 3^3\cdot 5\cdot 7\cdot 11\cdot 43$ & 86\\
$15561000=2^3\cdot 3^2\cdot 5^3\cdot 7\cdot 13\cdot 19$ & 76\\
$18763200=2^4\cdot 3^3\cdot 5^2\cdot 7\cdot 13\cdot 19$ & 76\\
$19061280=2^5\cdot 3^2\cdot 5\cdot 7\cdot 31\cdot 61$ & 61\\
$19845000=2^3\cdot 3^4\cdot 5^4\cdot 7^2$ & 40\\
$20638800=2^4\cdot 3^4\cdot 5^2\cdot 7^2\cdot 13$ & 48\\
$20884500=2^2\cdot 3^3\cdot 5^3\cdot 7\cdot 13\cdot 17$ & 85\\
$22932000=2^5\cdot  3^2\cdot 5^3\cdot 7^2\cdot 13$ & 64\\
$23152500^*=2^2\cdot 3^3\cdot 5^4\cdot 7^3$ & 49\\
$23569920=2^9\cdot 3^3\cdot 5\cdot 11\cdot 31$ & 80\\
$23647680=2^6\cdot 3^3\cdot 5\cdot 7\cdot 17\cdot 23$ & 92\\
$24160500=2^2\cdot 3^2\cdot 5^3\cdot 7\cdot 13\cdot 59$ & 59\\
$25798500=2^2\cdot 3^4\cdot 5^3\cdot 7^2\cdot 13$ & 54\\
$25832520=2^3\cdot 3^4\cdot 5\cdot 7\cdot 17\cdot 67$ & 67\\
$27027000=2^3\cdot 3^3\cdot 5^3\cdot 7\cdot 11\cdot 13$ & 110\\
$29381625=3^2\cdot 5^3\cdot 7^2\cdot 13\cdot 41$ & 41\\
$31872960=2^6\cdot 3^2\cdot 5\cdot 7\cdot 17\cdot 31$ & 93\\
$31888080=2^4\cdot 3^4\cdot 5\cdot 7\cdot 19\cdot 37$ & 74\\
$32997888=2^9\cdot 3^3\cdot 7\cdot 11\cdot 31$ & 84\\
$34889400=2^3\cdot 3^3\cdot 5^2\cdot 7\cdot 13\cdot 71$ & 71\\
$35626500=2^2\cdot 3^3\cdot 5^3\cdot 7\cdot 13\cdot 29$ & 87\\
$38383800=2^3\cdot 3\cdot 5^2\cdot 7\cdot 13\cdot 19$ & 74\\
$38785500=2^2\cdot 3^3\cdot 5^3\cdot 13^2\cdot 17$ & 52\\
$42997500=2^2\cdot 3^3\cdot 5^4\cdot 7^2\cdot 13$ & 52\\
$43205568=2^6\cdot 3\cdot 7\cdot 17\cdot 31\cdot 61$ & 61\\
$43330560=2^{10}\cdot 3\cdot 5\cdot 7\cdot 13\cdot 31$ & 80\\
$43857450=2\cdot 3^4\cdot 5^2\cdot 72\cdot 13\cdot 17$ & 51\\
$46683000=2^3\cdot 3^3\cdot 5^3\cdot 7\cdot 13\cdot 19$ & 114\\
$47297250=2\cdot 3^3\cdot 5^3\cdot 7^2\cdot 11\cdot 13$ & 77\\
$47392800=2^5\cdot 3\cdot 5^2\cdot 7^2\cdot 13\cdot 31$ & 62\\
$48323520=2^6\cdot 3^3\cdot 5\cdot 7\cdot 17\cdot 47$ & 94\\
$50213520=2^4\cdot 3^7\cdot 5\cdot 7\cdot 41$ & 72\\
$51597000=2^3\cdot 3^4\cdot 5^3\cdot 7^2\cdot 13$ & 72\\
$519795200=2^6\cdot 3\cdot 5^2\cdot 7^2\cdot 13\cdot 17$ & 64\\
$56511000=2^3\cdot 3^3\cdot 5^3\cdot 7\cdot 13\cdot 23$ & 115\\
$64701000=2^3\cdot 3^2\cdot 5^3\cdot 7\cdot 13\cdot 79$ &79\\
$68796000=2^5\cdot 3^3\cdot 5^3\cdot 7^2\cdot 13$ & 96\\
$71253000=2^3\cdot 3^3\cdot 5^3\cdot 7\cdot 13\cdot 29$ & 116\\
$77477400=2^3\cdot 3^2\cdot 5^2\cdot 7\cdot 11\cdot 13\cdot 43$ & 86\\
$77641200=2^4\cdot 3^3\cdot 5^2\cdot 7\cdot 13\cdot 79$ & 79\\
$93284100=2^2\cdot 3^2\cdot 5^2\cdot 7\cdot 13\cdot 17\cdot 67$ & 67\\
$95327232=2^{10}\cdot 3\cdot 7\cdot 11\cdot 13\cdot 31$ & 88\\
$98993664=2^9\cdot 3^4\cdot 7\cdot 11\cdot 31$ & 90\\
$103194000=2^4\cdot 3^4\cdot 5^3\cdot 7^2\cdot 13$ & 80\\
$108421632=2^9\cdot 3^3\cdot 11\cdot 23\cdot 31$ & 92\\
$109147500=2^2\cdot 3^4\cdot 5^4\cdot 7^2\cdot 11$ & 55\\
$109336500=2^2\cdot 3^3\cdot 5^3\cdot 7\cdot 13\cdot 89$ & 89\\
$129991680=2^{10}\cdot 3^2\cdot 5\cdot 7\cdot 13\cdot 31$ & 96\\
$133660800=2^7\cdot 3^3\cdot 5^2\cdot 7\cdot 13\cdot 17$ & 128\\
$136732050=2\cdot 3^4\cdot 5^2\cdot 7^2\cdot 13\cdot 53$ & 53\\
$138915000^*=2^3\cdot 3^4\cdot 5^4\cdot 7^3$ & 70\\
$142990842^9\cdot 3^2\cdot 7\cdot 11\cdot 13\cdot 31$ & 104\\
$144471600=2^4\cdot 3^4\cdot 5^2\cdot 7^3\cdot 13$ & 84\\
$144963000=2^3\cdot 3^3\cdot 5^3\cdot 7\cdot 13\cdot 59$ & 118\\
$160254000=2^5\cdot 3^2\cdot 5^3\cdot 7^3\cdot 13$ & 112\\
$164989440=2^9\cdot 3^3\cdot 5\cdot 7\cdot 11\cdot 31$ & 140\\
$172972800=2^8\cdot 3^3\cdot 5^2\cdot 7\cdot 11\cdot 13$ & 128\\
$176289750=2\cdot 3^3\cdot 5^3\cdot 7^2\cdot 13\cdot 41$ & 82\\
$188527500=2^2\cdot 3^4\cdot 5^4\cdot 7^2\cdot 19$ & 57\\
$191237760=2^7\cdot 3^4\cdot 5\cdot 7\cdot 17\cdot 31$ & 124\\
$199320576=2^{10}\cdot 3\cdot 7\cdot 13\cdot 23\cdot 31$ & 92\\
$219287250=2\cdot 3^4\cdot 5^3\cdot 7^2\cdot 13\cdot 17$ & 85\\
$221557248=2^9\cdot 3^3\cdot 11\cdot 31\cdot 47$ & 94\\
$227026800=2^4\cdot 3^4\cdot 5^2\cdot 7^2\cdot 11\cdot 13$ & 88\\
$232432200=2^3\cdot 3^3\cdot 5^2\cdot 7\cdot 11\cdot 13\cdot 43$ & 129\\
$247484160=2^8\cdot 3^4\cdot 5\cdot 7\cdot 11\cdot 31$ &124\\
$271498500=2^2\cdot 3^3\cdot 5^3\cdot 7\cdot 13^2\cdot 17$ & 91\\
$283783500=2^2\cdot 3^4\cdot 5^3\cdot 7^2\cdot 11\cdot 13$ & 99\\
$287752500=2^2\cdot 3^4\cdot 5^4\cdot 7^2\cdot 29$ & 58\\
$287878500=2^2\cdot 3^2\cdot 5^3\cdot 13\cdot 19\cdot 37$ & 111\\
$288943200=2^5\cdot 3^4\cdot 5^2\cdot 7^3\cdot 13$ & 108\\
$300982500=2^2\cdot 3^3\cdot 5^4\cdot 7^3\cdot 13$ & 91\\
$325798200=2^3\cdot 3^4\cdot 5^2\cdot 7\cdot 13^2\cdot 17$ & 78\\
$341775000=2^3\cdot 3^2\cdot 5^5\cdot 7^2\cdot 31$ & 70\\
$356879250=2\cdot 3^3\cdot 5^3\cdot 7^2\cdot 13\cdot 83$ & 83\\
$361179000=2^3\cdot 3^4\cdot 5^3\cdot 7^3\cdot 13$ & 126\\
$363854400=2^6\cdot 3\cdot 5^2\cdot 7^3\cdot 13\cdot 17$ & 112\\
$374078250=2\cdot 3^4\cdot 5^3\cdot 7^2\cdot 13\cdot 29$ & 87\\
$377055000=2^3\cdot 3^4\cdot 5^4\cdot 7^2\cdot 19$ & 76\\
$389975040=2^{10}\cdot 3^3\cdot 5\cdot 7\cdot 13\cdot 31$ & 144\\
$390957840=2^4\cdot 3^5\cdot 5\cdot  7\cdot  13^2\cdot 17$ & 104\\
$407307264=2^{10}\cdot 3\cdot 7\cdot 13\cdot 31\cdot 47$ & 94\\
$421866900=2^2\cdot 3^3\cdot 5^2\cdot 7\cdot 13\cdot 17\cdot 101$ & 101\\
$428972544=2^9\cdot 3^3\cdot 7\cdot 11\cdot 13\cdot 31$ & 156\\
$434397600=2^5\cdot 3^3\cdot 5^2\cdot 7\cdot 13^2\cdot 17$ & 104\\
$438574500=2^2\cdot 3^4\cdot 5^3\cdot 7^2\cdot 13\cdot 17$ & 102\\
$447828480=2^9\cdot 3^3\cdot 5\cdot 11\cdot 19\cdot 31$ & 152\\
$467002900=2^2\cdot 5^2\cdot 7\cdot 13\cdot 19\cdot 37\cdot 73$ & 73\\
$474692400=2^4\cdot 3^4\cdot 5^2\cdot 7^2\cdot 13\cdot 23$ & 92\\
$481572000=2^5\cdot 3^3\cdot 5^3\cdot 7^3\cdot 13$ & 168\\
$486319680=2^6\cdot 3^3\cdot 5\cdot 7\cdot 11\cdot 17\cdot 43$ & 172\\
$488697300=2^2\cdot 3^5\cdot 5^2\cdot 7\cdot 13^2\cdot 17$ & 81\\
$490990500=2^2\cdot 3^2\cdot 5^3\cdot 7\cdot 11\cdot 13\cdot 109$ & 109\\
$494968320=2^9\cdot 3^4\cdot 5\cdot 7\cdot 11\cdot 31$ & 150\\
$513513000=2^3\cdot 3^3\cdot 5^3\cdot 7\cdot 11\cdot 13\cdot 19$ & 209\\
$552348720=2^4\cdot 3^7\cdot 5\cdot 7\cdot 11\cdot 41$ & 132\\
$555660000^*=2^5\cdot 3^4\cdot 5^4\cdot 7^3$ & 100\\
$559704600=2^3\cdot 3^3\cdot 5^2\cdot 7\cdot 13\cdot 17\cdot 67$ & 134\\
$567567000=2^3\cdot 3^4\cdot 5^3\cdot 7^2\cdot 11\cdot 13$ & 132\\
$575757000=2^3\cdot 3^2\cdot 5^3\cdot 7\cdot 13\cdot 19\cdot 37$ & 148\\
$585427500=2^2\cdot 3^4\cdot 5^4\cdot 7^2\cdot 59$ & 59\\
$639802800=2^4\cdot 3^4\cdot 5^2\cdot 7^2\cdot 13\cdot 31$ & 93\\
$648083520=2^6\cdot 3^2\cdot 5\cdot 7\cdot 17\cdot 31\cdot 61$ & 122\\
$648784500=2^2\cdot 3\cdot 5^3\cdot 7^3\cdot 13\cdot 97$ & 97\\
$690908400=2^4\cdot 3^3\cdot 5^2\cdot 7\cdot 13\cdot 19\cdot 37$ & 148\\
$708107400=2^3\cdot 3^3\cdot 5^2\cdot 7\cdot 11\cdot 13\cdot 131$ & 131\\
$710892000=2^5\cdot 3^2\cdot 5^3\cdot 7^2\cdot 13\cdot 31$ & 124\\
$722358000=2^4\cdot 3^4\cdot 5^3\cdot 7^3\cdot 13$ & 140\\
$756756000=2^5\cdot 3^3\cdot 5^3\cdot 7^2\cdot 11\cdot 13$ & 176\\
$758951424=2^9\cdot 3^3\cdot 7\cdot 11\cdot 23\cdot 31$ & 161\\
$779688000=2^6\cdot 3^2\cdot 5^3\cdot 7^2\cdot 13\cdot 17$ & 128\\
$793457920=2^7\cdot 3^4\cdot 5\cdot 7\cdot 17\cdot 127$ & 127\\
$823280640=2^{10}\cdot 3\cdot 5\cdot 7\cdot 13\cdot 19\cdot 31$ & 152\\
$953629840=2^4\cdot 3^7\cdot 5\cdot 7\cdot 17\cdot 41$ & 136\\
$877149000=2^3\cdot 3^4\cdot 5^3\cdot 7^2\cdot 13\cdot 17$ & 136\\
$879196500=2^2\cdot 3^2\cdot 5^3\cdot 7\cdot 13\cdot 19\cdot 113$ & 113\\
$970023600=2^4\cdot 3^4\cdot 5^2\cdot 7^2\cdot 13\cdot 47$ & 94\\
$973176750=2\cdot 3^2\cdot 5^3\cdot 7^3\cdot 13\cdot 97$ & 97\\
$977394600=2^3\cdot 3^5\cdot 5^2\cdot 7\cdot 13^2\cdot 17$ & 108\\
$992548080=2^4\cdot 3^8\cdot 5\cdot 31\cdot 61$ & 81\\ \hline
\end{longtable}
\ec


\begin{thebibliography}{99}
\bibitem{1}
G. L. Cohen,
{\it Numbers whose positive divisors have small integral harmonic mean},
Math. Comp. {\bf 66}(1997), 883-891.

\bibitem{2}
G. L. Cohen and R. M. Sorli,
{\it Harmonic seeds},
Fib. Quart. {\bf 36}(1998), 386-390.

\bibitem{3}
G. L. Cohen and M. Deng,
{\it On a generalization of Ore's harmonic numbers},
Nieuw Arch. Wiskunde {\bf 16}(1998), no. 3, 161-172.

\bibitem{4}
M. Garcia,
{\it On numbers with integral harmonic mean},
Amer. Math. Monthly, {\bf 61}(1954), 89-96.

\bibitem{18}
T. Goto and S. Shibata,
{\it All numbers whose positive divisors have integral harmonic mean
up to 300},
Math. Comp. {\bf 73}(2004), 475-491.

\bibitem{5}
P. Hagis and G. Lord,
{\it Unitary harmonic numbers},
Proc. Amer. Math. Soc. {\bf 51}(1975), 1-7.

\bibitem{6}
P. Hagis and G. L. Cohen,
{\it Infinitary harmonic numbers},
Bull. Austral. Math. Soc. {\bf 41}(1990), 151-158.

\bibitem{17}
H.-J. Kanold,
{\it \"Uber das harmonische Mittel der Teiler einer nat\"urlichen
Zahl},
Math. Ann. {\bf 13}(1957), 371-374.

\bibitem{7}
W. H. Mills,
{\it On a conjecture of Ore},
Proc. Number Theory Conf., Boulder Co., 1972, 142-146.

\bibitem{8}
K. Nageswara Rao,
{\it On some unitary divisor functions},
Scripta Math. {\bf 28}(1967), 347-351.

\bibitem{9}
O. Ore,
{\it On the averages of the divisors of a number},
Amer. Math. Monthly, {\bf 55}(1948), 615-619.

\bibitem{10}
C. Pomerance,
{\it On a problem of Ore: Harmonic numbers},
Abstract 709-A5, Notices Amer. Math. Soc. {\bf 20}(1973), A-648.

\bibitem{19}
J. S\'andor,
{\it Handbook of number theory}, II, (In coop. with B. Crstici),
Springer Verlag, 2004.

\bibitem{11}
M. V. Subbarao,
{\it Problem E1558},
Amer. Math. Monthly {\bf 70}(1963), 92, solution in
{\bf 70}(1963), 1009-1010.

\bibitem{12}
M. V. Subbarao and L. J. Warren,
{\it Unitary perfect numbers},
Canad. Math. Bull. {\bf 9}(1966), 147-153.

\bibitem{13}
D. Suryanarayana,
{\it The number of bi-unitary divisors of an integer},
Lecture Notes in Math., vol. 251, 1972, 273-278.

\bibitem{14}
Ch. Wall,
{\it Bi-unitary perfect numbers},
Proc. Amer. Math. Soc. {\bf 33}(1972), no. 1, 39-42.

\bibitem{15}
Ch. R. Wall,
{\it The fifth unitary perfect number},
Canad. Math. Bull. {\bf 18}(1975), 115-122.

\bibitem{16}
Ch. Wall,
{\it Unitary harmonic numbers},
Fib. Quart. {\bf 21}(1983), 18-25.






\end{thebibliography}
\end{document}